\documentclass{article}
\usepackage{amsfonts}
\usepackage{amsmath}

\newtheorem{theorem}{Theorem}

\newtheorem{conclusion}[theorem]{Conclusion}

\newtheorem{corollary}[theorem]{Corollary}

\newtheorem{proposition}[theorem]{Proposition}
\newtheorem{remark}[theorem]{Remark}

\newenvironment{proof}[1][Proof]{\noindent\textbf{#1.} }{\ \rule{0.5em}{0.5em}}
\input{tcilatex}

\begin{document}

\title{Global Existence of solutions for systems of coupled reaction
diffusion equations with nonlinearities of unlimited growth.}
\author{Said Kouachi \\
Abbes Laghrour University, Khenchela, Algeria.}
\maketitle

\begin{abstract}
In this work we prove global existence and uniform boundedness of solutions
of 2$\times $2 reaction-diffusion systems with control of mass structure and
nonlinearities of unlimited growth. Furthermore the results are obtained
without restrictions on diffusion terms neither on the initial data. Such
systems possess many applications in physical-chemistry. Our technique of
proof relies on a judiciously rectified Lyapunov functional used previously
by the author in several papers
\end{abstract}

\section{\protect\bigskip \textbf{Introduction}}

In this article, we are concerned with the existence of globally bounded
solutions to the reaction-diffusion system 
\begin{equation}
\left\{ 
\begin{array}{c}
\\ 
\frac{\partial u}{\partial t}-a\Delta u=f(u,v), \\ 
\\ 
\frac{\partial v}{\partial t}-b\Delta v=g(u,v), \\ 
\end{array}%
\right. \quad \mathrm{in}\quad \mathbb{R}^{+}\times \Omega ,  \label{1.1}
\end{equation}%
subject to the boundary conditions 
\begin{equation}
\frac{\partial u}{\partial \eta }=\frac{\partial v}{\partial \eta }=0\qquad
on\quad \mathbb{R}^{+}\times \partial \Omega  \label{1.2}
\end{equation}%
and the initial data 
\begin{equation}
u(0,x)=u_{0}(x),\qquad v(0,x)=v_{0}(x)\qquad \mathrm{in}\quad \Omega
\label{1.3}
\end{equation}%
where $a>0$ and $b>0$ are the diffusion coefficients of some interacting
species whose spatiotemporal densities are $u$ and $v$. The domain $\Omega $
is an open bounded domain of class $C^{1}$ in $\mathbb{R}^{n}$, with
boundary $\partial \Omega $ and $\frac{\partial }{\partial \eta }$ denotes
the outward derivative on $\partial \Omega $. The initial data are assumed
to be nonnegative and bounded. The reactions $f$ \ and $g$ are continuously
differentiable functions. We suppose that, for some positive constants $C$
and $\mu $ we have%
\begin{equation}
f(u,v)\leq f(u,v)+\mu g(u,v)\leq 0,\text{ for all }u\geq 0,\ v\geq 0\text{
with }u+v\geq C.  \label{1.5}
\end{equation}%
The last inequality is called "the control of mass condition".\newline
We suppose $f(0,v)=0$ and $g(u,0)\geq 0,$ for all $u\geq 0$ and$\ v\geq 0$,
then using standard comparison arguments for parabolic equations (see e.g. 
\cite{Hen}) the solutions remain positive at any time$.$

In the case of systems on the form%
\begin{equation}
\left\{ 
\begin{array}{c}
\\ 
\frac{\partial u}{\partial t}-a\Delta u=-uF(v), \\ 
\\ 
\frac{\partial v}{\partial t}-b\Delta v=uG(v), \\ 
\end{array}%
\right. \quad \mathrm{in}\quad \mathbb{R}^{+}\times \Omega ,  \label{1.6}
\end{equation}%
\cite{Har-You} generalized the results of K. Masuda for $F(v)=G(v)=v^{\beta
} $ to nonlinearities $F(s)$ satisfying the condition 
\begin{equation}
\displaystyle\lim_{s\rightarrow +\infty }\left[ \frac{\ln \left(
1+F(s)\right) }{s}\right] =0,  \label{1.7}
\end{equation}%
which extends the nonlinearities to sub-linear ones (e.g. $%
F(v)=G(v)=e^{v^{\gamma }},\;0<\gamma <1)$. In the case of nonlinearities of
exponential growth like $F(v)=G(v)=e^{v}$ which appears in the
Frank-Kamenetskii approximation to Arrhenius-type reaction \cite{Ari}, \cite%
{Bar} has made a small progress in this direction and proved that solutions
are globally bounded under the condition $||u_{0}||_{\infty }\leq
8ab/(a-b)^{2}$. Then \cite{Mar-Pie} proved (in the case $0<b<a<\infty $
which means that the absorbed substance diffuses faster than the other one)
that global solutions exist if $\Omega =\mathbb{R}^{n}$. The proof is based
on a simple comparison property concerning the kernels associated with the
operators $\left( \frac{\partial }{\partial t}-a\Delta \right) $ and $\left( 
\frac{\partial }{\partial t}-b\Delta \right) .$ Some years ago, \cite%
{Her-Lac-Vel} obtained similar results (in the general case). Then the
authors in \cite{Kou-You} generalized the results obtained in \cite{Har-You}
while adding $-c\Delta u$ to the right-hand side of the second equation of
system (\ref{1.1}) under the condition \ 
\begin{equation*}
\underset{s\rightarrow +\infty }{lim}\left[ \frac{log\left( 1+f(r,s)\right) 
}{s}\right] <\alpha ^{\ast }\text{, for\ any}\;r\geq 0,
\end{equation*}%
where $f(r,s)\geq 0$ for all $r,\ s\geq 0$ and

\begin{equation*}
\alpha ^{\ast }=\frac{2ab}{n(a-b)^{2}\left\Vert u_{0}\right\Vert _{\infty }},
\end{equation*}%
condition reflecting the weak exponential growth of the reaction term $f$. 
\cite{Kir-Kan 3} and then \cite{Hos} proved global existence of solutions of
(\ref{1.6}) on bounded domains if $a>b$ or $c\geq a-b>0$, while in \cite%
{Kir-Kan 4} it is showed that solutions exist\ globally whenever $b>a,\
c<b-a $ and $G(v)$ is of exponential growth. Recently in \cite{Kou(MMAS)} we
proved global existence of solutions of (\ref{1.6}) in a bounded domain
under the condition

\begin{equation}
\displaystyle\lim_{v\rightarrow +\infty }\frac{G^{\prime }(v)}{F(v)}=0,
\label{1.8}
\end{equation}%
where $G^{\prime }$ denotes the first derivative of $G$ with respect to $v$
and the functions $F$, $G$, $G^{\prime }$ and $G^{\prime \prime }$are
continuously differentiable and nonnegative. \newline
To the best of our knowledge, the question of global existence of reaction
diffusion systems on a bounded domain remains open when the reactions grow
faster than a polynomial. Some partial positive results have been obtained
only when the reactions grow faster than a polynomial as it is cited above
(see \cite{Kou(MMAS)} and \cite{Bar} when $\Omega $ is bounded and \cite%
{Mar-Pie} in the unbounded case or when $\Omega =%
\mathbb{R}
^{n}$). We should know that, in some cases, solutions of system (\ref{1.1})
with appropriate boundary conditions can plow-up in finite time (see for
example \cite[ 22]{Pie-Sch 22}).

In this paper we show the global existence of a unique solution (uniformly
bounded on $\mathbb{R}^{+}\times\Omega$) to problem (\ref{1.1})-(\ref{1.3})
without conditions on the nonlinearities growth (i.e. condition (\ref{1.5}))
as it was supposed until now. Remark that for system (\ref{1.6}), conditions
like (\ref{1.8}) are more stronger then (\ref{1.5}) when the reactions grow
faster than a polynomial, especially when $F=G.$

\section{\textbf{\protect\bigskip Statement of the main results}}

It is well known that, for any initial data in $L^{\infty }(\Omega )$, local
existence and uniqueness of classical solutions to the initial value problem
(\ref{1.1})-(\ref{1.3}) follows from the basic existence theory for abstract
semi-linear differential equations (see \cite{Ama 85}, \cite{Hen} or \cite%
{Rot}). The solutions are classical on $(0,T_{\max }),$ where $T_{\max }$
denotes the eventual blowing-up time in $L^{\infty }(\Omega )$. If we denote
by%
\begin{equation*}
\left\Vert u\right\Vert _{\infty }\;=\;\max_{x\in \Omega }\left\vert
u(x)\right\vert ,
\end{equation*}%
the usual norm in the space $L^{\infty }(\Omega )\ $(or $C\left( \overline{%
\Omega }\right) $), let us recall the following classical local existence
result under the above assumptions:

\begin{proposition}
\label{2}The problem(\ref{1.1})-(\ref{1.3}) admits a unique, classical
solution $(u,v)$ on $(0,T_{\max }[\times \Omega $. If $T_{\max }<\infty $
then 
\begin{equation*}
\lim_{t\nearrow T_{\max }}\left\{ \left\Vert u(t,.)\right\Vert _{\infty
}+\left\Vert v(t,.)\right\Vert _{\infty }\right\} =\infty .
\end{equation*}
\end{proposition}

Then we can formulate our main result as follows:

\begin{theorem}
\label{1}Under condition (\ref{1.5}), the solution of problem (\ref{1.1})-(%
\ref{1.3}) with positive initial data in $\mathbb{L}^{\infty }(\Omega )$
exists globally in time. Moreover $u$ and $v$ are uniformly bounded with%
\begin{equation}
\left\Vert u(t,.)\right\Vert _{\infty }\leq \overline{u}_{0}\text{ and \ }%
\left\Vert v(t,.)\right\Vert _{\infty }\leq \overline{v}_{0},\text{ for all }%
t>0,  \label{2.1}
\end{equation}%
where%
\begin{equation}
\overline{u}_{0}=\sup \left\{ C,\left\Vert u_{0}\right\Vert _{\infty
}\right\} \text{ and }\overline{v}_{0}=\sup \left\{ C,\left\Vert
v_{0}\right\Vert _{\infty }\right\} ,  \label{2.2}
\end{equation}%
and $C$ given by (\ref{1.5}).
\end{theorem}

For example our results are applicable to system (\ref{1.6}) for $%
F(v)=e^{e^{v}}-P(v)$ and $G(v)=e^{e^{v}},$ where $P$ \ is any polynomial
(even zero polynomial) of finite degree$.$\newline

This result can be formulated by the following

\begin{corollary}
\label{5}The solutions of system (\ref{1.6}) with boundary conditions (\ref%
{1.2}) and positive initial data in $\mathbb{L}^{\infty }(\Omega )$ are
global and uniformly bounded on $\left[ 0,+\infty \right[ \times \Omega $
with%
\begin{equation*}
\left\Vert u(t,.)\right\Vert _{\infty }\leq \left\Vert u_{0}\right\Vert
_{\infty }\text{ and \ }\left\Vert v(t,.)\right\Vert _{\infty }\leq 
\overline{v}_{0},\text{ for all }t>0,
\end{equation*}%
provided that%
\begin{equation*}
0<\lim_{s\rightarrow +\infty }\frac{F(s)}{G(s)}\leq +\infty .
\end{equation*}
\end{corollary}

Another application of Theorem \ref{1}, is the following system modelling
the exothermic combustion in a gas%
\begin{equation}
\left\{ 
\begin{array}{ccc}
\frac{\partial Y}{\partial t}-a\Delta Y & = & -Y^{m}e^{T}, \\ 
\frac{\partial T}{\partial t}-b\Delta T & = & Y^{m}e^{T},%
\end{array}%
\right. \quad \mathrm{in}\quad \mathbb{R}^{+}\times \Omega ,  \label{2.5}
\end{equation}%
where $Y$ \ is the concentration of a single reactant, $T$ \ is the
temperature (see e.g. \cite{Pie} in the case when the function $e^{T}$ is
approximated by a polynomial and \cite{Her-Lac-Vel} in more general case).
Moreover, we have the following result

\begin{corollary}
\label{6}The system (\ref{2.5}) with boundary conditions (\ref{1.2}) and
positive initial data in $\mathbb{L}^{\infty }(\Omega )$ has global and
uniformly bounded solutions on $\left[ 0,+\infty \right[ \times \Omega $
satisfying%
\begin{equation*}
\left\Vert Y(t,.)\right\Vert _{\infty }\leq \left\Vert Y(0,.)\right\Vert
_{\infty }\text{ and \ }\left\Vert T(t,.)\right\Vert _{\infty }\leq
\left\Vert T(0,.)\right\Vert _{\infty },\text{ for all }t\geq 0.
\end{equation*}
\end{corollary}

\section{The \textbf{Proofs}}

For the proofs, we use a judiciously rectified Lyapunov functional used in 
\cite{Kou(ED2001)}, \cite{Par-Kou-Gut}, \cite{Kou Arx 22} and \cite%
{Kou(DPDE)}.

Let $\theta >1$ be any positive constant satisfying%
\begin{equation}
\theta ^{2}>\frac{\left( a+b\right) ^{2}}{4ab},  \label{3.1}
\end{equation}%
and let $\theta _{0},\ \theta _{1},...,\theta _{p}$ a real positive numbers
satisfying%
\begin{equation}
\dfrac{\theta _{i}\theta _{i+2}}{\theta _{i+1}^{2}}=\theta ^{2},\ \ \ \
i=0,1,...,p.  \label{3.2}
\end{equation}%
where $p$ is a positive integer. If we denote by $w_{+}$ the nonnegative
part $sup\{w,0\}$ of a function $w$, then define the functional 
\begin{equation}
t\longrightarrow L(t)=\int\limits_{\Omega }H\left( u(t,x),v(t,x)\right) dx,
\label{3.3}
\end{equation}%
\newline
where%
\begin{equation}
H(u,v)=\overset{p}{\underset{i=0}{\sum \ }}C_{p}^{i}\theta _{i}U^{i}V^{p-i},
\label{3.4}
\end{equation}%
\newline
where%
\begin{equation}
U=\left( u-\overline{u}_{0}\right) _{+}\text{ and }V=\left( v-\overline{v}%
_{0}\right) _{+},  \label{3.5}
\end{equation}%
with%
\begin{equation*}
C_{p}^{i}=\frac{p!}{i!(p-i)!},\ i=0,1,...,p,
\end{equation*}%
and $i!$ denotes the product of positive integers less then $i$.

\begin{proof}
\textbf{(of Theorem \ref{1})}By differentiating $L$ with respect to $t$ and
then by simple use of Green's formula we get%
\begin{equation*}
L^{\prime }(t)=I+J,
\end{equation*}%
where 
\begin{equation}
I=-p(p-1)\overset{p-2}{\underset{i=0}{\tsum }}C_{p-2}^{i}\
\dint\limits_{\Omega }T_{i}\left( \nabla u,\nabla v\right) U^{i}V^{p-2-i}dx,
\label{3.6}
\end{equation}%
with for all $i=0,...,p-1$%
\begin{equation*}
T_{i}\left( \nabla u,\nabla v\right) =\left( a\theta _{i+2}sgn\left(
U\right) \left\vert \nabla u\right\vert ^{2}+\left( a+b\right) \theta
_{i+1}sgn\left( U\right) sgn\left( V\right) \nabla u\nabla v+b\theta
_{i}sgn\left( V\right) \left\vert \nabla v\right\vert ^{2}\right) ,
\end{equation*}%
and 
\begin{equation}
J=p\overset{p-1}{\underset{i=0}{\sum }}C_{p-1}^{i}\int\limits_{\Omega
}\left( \theta _{i+1}sgn\left( U\right) f(u,v)+\theta _{i}sgn\left( V\right)
g(u,v)\right) U^{i}V^{p-1-i}dx.\newline
\label{3.7}
\end{equation}%
Then condition (\ref{3.2}) gives the positivity of the quadratics $%
T_{i}\left( \nabla u,\nabla v\right) $ to get $I\leq 0$ on $\left( 0,T_{\max
}\right) $.\newline
Let us choose $\frac{\theta _{0}}{\theta 1}<\mu ,\ $then from (\ref{3.1})-(%
\ref{3.2}) $\frac{\theta _{i}}{\theta _{i+1}}<\mu $ ; for all $i=0,...,p-1$.
We consider the following two alternatives:\newline
(i) If%
\begin{equation}
u+v\geq C,  \label{3.8}
\end{equation}%
then, since in this case and from condition (1.4) $g$ is nonnegative, we have%
\begin{equation*}
f(u,v)+\mu g(u,v)\leq 0,\text{ for all }(u,v)\text{ satisfying (\ref{3.8}).}
\end{equation*}%
This gives%
\begin{equation*}
\left( \theta _{i+1}f(u,v)+\theta _{i}g(u,v)\right) U^{i}V^{p-1-i}\leq 0;%
\text{ for all\ }i=0,...,p-1.
\end{equation*}%
(ii) If%
\begin{equation*}
u+v<C,
\end{equation*}%
and since $u$ and $v$ are nonnegative, we get from (\ref{2.2})%
\begin{equation*}
u<\overline{u}_{0}\text{ and }v<\overline{v}_{0},
\end{equation*}%
which gives%
\begin{equation*}
U=\left( u-\overline{u}_{0}\right) _{+}\text{ }=\text{ }V=\left( v-\overline{%
v}_{0}\right) _{+}=0,
\end{equation*}%
and consequently%
\begin{equation*}
\left( \theta _{i+1}f(u,v)+\theta _{i}g(u,v)\right) U^{i}V^{p-1-i}=0;\text{
for all\ }i=0,...,p-1.
\end{equation*}%
This completes this alternative and gives $J\leq 0$ on $\left( 0,T_{\max
}\right) $. So the nonnegative functional $L(t)$ is decreasing on $\left(
0,T_{\max }\right) $\ and since $L(0)=0$, this gives%
\begin{equation*}
L(t)=0\text{, for all }t<T_{\max }.
\end{equation*}%
Thus 
\begin{equation*}
U=V=0,\text{ for all }x\in \Omega \text{ and all }t<T_{\max },
\end{equation*}%
and then%
\begin{equation*}
u\leq \overline{u}_{0}\text{ and }v\leq \overline{v}_{0},\text{ for all }%
x\in \Omega \text{ and all }t<T_{\max }.
\end{equation*}%
The application of Proposition \ref{2} gives $T_{\max }=+\infty $ and ends
the proof of the Theorem.\newline
\end{proof}

\begin{proof}
\textbf{(of Corollary \ref{5})} Since the reaction corresponding to the
first equation of (\ref{1.6}) is nonpositive, then by direct application of
the maximum principle we have%
\begin{equation*}
u\left( t,x\right) \leq \left\Vert u_{0}\right\Vert _{\infty },\text{ for
all }x\in \Omega \text{ and all }t<T_{\max }.
\end{equation*}%
If we take any positive constant $\lambda <1$ such that%
\begin{equation*}
\lim_{s\rightarrow +\infty }\frac{F(s)}{G(s)}>\lambda ,
\end{equation*}%
then we can find a positive constant $A$ such that%
\begin{equation*}
\frac{F(s)}{G(s)}>\lambda ,\text{ for all }s>A.
\end{equation*}%
The Corollary is then an immediate consequence of the Theorem \ref{1} by
taking%
\begin{equation*}
C=A.
\end{equation*}
\end{proof}

\begin{proof}
\textbf{(of Corollary \ref{6})} The proof is similar to that of Corollary %
\ref{5} by taking $C=0$ and any positive constant $\lambda <1$ .
\end{proof}

\begin{remark}
Theorem \ref{1} is applicable for more general boundary conditions than (\ref%
{1.2}) as homogenous and nonhomogenous Dirichlet and mixed boundary
conditions.
\end{remark}

\begin{remark}
Condition (\ref{1.5}) should be taken with prudence: It can arrive that it
is not satisfied for systems with simple reaction of polynomial growth. It
is the case, for example when $f(u,v)=\left( u-u^{2}\right) v^{2}$ and $%
g(u,v)=uv^{2}$ where the solution blows up at finite time for initial data $%
\frac{1}{2}\leq u_{0}\leq 1$ on $\Omega .$ Indeed the region%
\begin{equation*}
\sum =\left\{ \left( u,v\right) \in 
\mathbb{R}
_{+}^{2}\text{ such that }\frac{1}{2}\leq u\leq 1\right\} ,
\end{equation*}%
is an inariant region which gives $g(u,v)\geq \frac{1}{2}v^{2}$ in $\sum $,
then by comparison arguments we get\ blow up at finite time. Obviously we
can not find positive constants $C$ and $\lambda $ such that (\ref{1.5}) is
satisfied.
\end{remark}

\begin{conclusion}
The global existence of solutions to systems of a coupled reaction diffusion
equations with unlimited nonlineraties growth had been open until very
recently except for some partial results under very restrictive conditions
as small initial data. The main reason to obtain interesting results in such
situation is the lack of the classical methods as those based on the entropy
inequality or those using the duality argument etc... In this work we used a
very simple functional but very strong solving the problem of global
existence of reaction diffusion systems without conditions on the
nonlineraties growth.
\end{conclusion}

\end{document}